\documentstyle[12pt]{article}
\def \Z{\hbox{$Z\hskip -5.2pt Z$}}
\def\OVER#1#2{\mbox{${{\ssc\,}\dis#1{\ssc\,}}
\over{{\ssc\,}\rb{-2pt}{\mbox{$\dis#2$}}{\ssc\,}}$}}
\def\LL{S_\mu}

\def \C{\hbox{$C\hskip -5pt \vrule height 6pt depth 0pt \hskip 6pt$}}
\def \sC{{\hbox{$\sc C\hskip -5pt \vrule height 5pt depth 0pt \hskip 6pt$}}}
\def\qed{\hfill \hfill \ifhmode\unskip\nobreak\fi\ifmmode\ifinner
         \else\hskip5pt\fi\fi
 \hbox{\hskip5pt\vrule width4pt height6pt depth1.5pt\hskip 1 pt}}
\def\a{\alpha}

\def\d{\delta}

\def\l{\lambda}

\def\SA{\hbox{$S\!A$}}
\def\SB{\hbox{$S\!B$}}
\def\Vir{\hbox{\bf\sl Vir}}
\def\HVir{{}}

\def\SVio{\hbox{{\bf\sl SVir}${\ssc\,}_0$}}
\def\SVii{\hbox{{\bf\sl SVir}${\ssc\,}_1$}}
\def\HVir{\hbox{{\bf\sl SVir}[$M$]}}
\def\SVir{\hbox{{\bf\sl SVir}[$M,s$]}}
\def\SVira{\hbox{{\bf\sl SVir}[$M_1,s$]}}

\def\sc{\scriptstyle}
\def\ssc{\scriptscriptstyle}
\def\dis{\displaystyle}
\def\cl{\centerline}

\def\nl{\newline}
\def\ol{\overline}

\def\rar{\rightarrow}

\def\Lra{\Leftrightarrow}
\def\bs{\backslash}
\def\hs{\hspace*}
\def\vs{\vspace*}

\def\rb{\raisebox}

\def\ni{\noindent}
\def\hi{\hangindent}
\def\ha{\hangafter}

\begin{document}
\def\TIT{SIMPLE MODULES OVER
THE HIGHER RANK SUPER-VIRASORO ALGEBRAS}
\def\ABS
{It is proved that uniformly bounded simple modules over higher rank
super-Virasoro algebras are modules of the intermediate series, and
that simple modules with finite dimensional weight spaces are either
modules of the intermediate series or generalized highest weight modules.}
\def\KEYW{higher rank super-Virasoro algebra,
uniformly bounded module,
module of the intermediate series, generalized highest weight module.}
\def\cover{{\bf\par\ \par\ \par\ni\hi5ex\ha1
Article title:\nl \TIT
\par\ \par\ \par\ni\hi5ex\ha1
Mailing address:\nl
  Dr. Yucai Su \nl
  Department of Applied Mathematics \nl
  Shanghai Jiaotong University \nl
  1954 Huashan Road, Shanghai 200030 \nl
  China \nl
  Email: kfimmi@public1.sta.net.cn
\par\ \par\ \par\ni\hi5ex\ha1
Running title:\nl
  Modules over higher rank super-Virasoro algebras
\par\ \par\ \par\ni\hi5ex\ha1
Abstract:\nl \ABS
\par\ \par\ni\hi5ex\ha1 Keywords:\nl
\hs{5ex}\KEYW
}\setcounter{page}{0}\pagebreak}
\cl{\bf SIMPLE MODULES OVER THE HIGHER RANK}
\cl{\bf SUPER-VIRASORO ALGEBRAS}
{\small\vs{4pt}\cl{(Appeared in {\it Lett.Math.Phys.} {\bf53}
(2000), 263-272)}
\par\vs{-7pt}\ \par
\cl{Yucai Su\footnote{Supported by a grant from National
Educational Department of China}} \vs{-12pt}
\par\ \par
\cl{\small\it Department of Mathematics, Shanghai Jiaotong
University, Shanghai 200030, China} \cl{ Email: ycsu@sjtu.edu.cn}
\vs{-7pt} \ \par \ni{\bf Abstract.} \ABS
\par\ni
\par\ni{\bf Keywords:} \KEYW
}\vs{-7pt}
\par\
\par
\cl{\bf \S1. Introduction}\par
The Virasoro algebra $\Vir$, playing a fundamental role in two dimensional
conformal quantum field theory, is the universal central extension of the
complex Lie algebra of the polynomial vector fields on the circle [1].
The notion of higher rank Virasoro algebras and the notion of higher rank
super-Virasoro algebras were introduced in [2,3] and the notions
of generalized Virasoro and super-Virasoro algebras were introduced in [4].
Kac [5] conjectured a theorem that
a simple $\Vir$-module with finite dimensional weight spaces is either
a module of the intermediate series or else a highest or lowest weight
module. This theorem was partially proved in [6,7,8] and fully proved
in [9] and then generalized to the super-Virasoro algebras in [10], and
was in some sense generalized to the higher rank Virasoro algebras in [11].
In this paper, we shall obtain a similar result for
the higher rank super-Virasoro algebras, thus completing the proof of
the theorem conjectured by the author in [3].
In [12,] by using the above mentioned theorem, we gave a theorem on
indecomposable $\Vir$-modules.
%
%
Using our result here, we shall have a better understanding on
the simple modules over the higher rank super-Virasoro algebras and
it may be possible to give in some extent a classification of the
uniformly bounded modules over the higher rank super-Virasoro algebras.
This is also our motivation to present the result here.
\par
Let $n$ be a positive integer.
Let $M$ be an $n$-dimensional $\Z$-submodule of $\C$, and let
$s\in\C$ such that $2s\in M$. The {\it rank $n$
super-Virasoro algebra} (or a {\it higher rank super-Virasoro algebra}
if $n\ge2$) is the Lie superalgebra $\SVir=\SVio\oplus\SVii$,
where $\SVio$ has a basis $\{L_\mu,c\,|\,\mu\in M\}$ and
$\SVii$ has a basis $\{G_\eta\,|\,\eta\in s+M\}$, with the commutation
relations
$$
\matrix{
[L_\mu,L_\nu]=(\nu-\mu)L_{\mu+\nu}-\d_{\mu+\nu,0}\OVER1{12}(\mu^3-\mu)c,
\hfill&
[L_\mu,G_\eta]=(\eta-\OVER{\mu}2)G_{\mu+\eta},
\vs{4pt}\hfill\cr
[G_\eta,G_\l]=2L_{\eta+\l}-\d_{\eta+\l,0}\OVER13(\eta^2-\OVER14)c,
\hfill&
[L_\mu,c]=[G_\eta,c]=0,
\hfill\cr}
\eqno(1.1)$$
for $\mu,\nu\in M,\,\eta,\l\in s+M.$
Then $\SVio$ is the higher rank Virasoro algebra $\HVir$.
It is proved in [4] that
a {\it module of the intermediate series} (i.e., a module with all the
dimensions of the weight spaces of ``even'' or ``odd'' part
are $\le1$) over $\SVir$ is one of the three
series of the modules $\SA_{a,b},\SA_{a'},\SB_{a'}$ or their quotient modules
for suitable $a,b,a'\in\C$, where $\SA_{a,b},\SA_{a'}$ have basis
$\{x_\mu\,|\,\mu\in M\}\cup\{y_\eta\,|\,\eta\in s+M\}$ and $\SB_{a'}$
has basis $\{x_\eta\,|\,\eta\in s+M\}\cup\{y_\mu\,|\,\mu\in M\}$ such that
$c$ acts trivially and such
that
$$
\matrix{
\matrix{
\SA_{a,b}\!:\!\!\!\!\hfill&
L_\mu x_\nu=(a+\nu+\mu b)x_{\mu+\nu},\hfill&
L_\mu y_\eta=(a+\eta+\mu(b-\OVER12))y_{\mu+\eta}
\vs{4pt}\hfill\cr&
G_\l x_\nu=y_{\l+\nu},\hfill&
G_\l y_\eta=(a+\eta+2\l(b-\OVER12))x_{\l+\eta}
\hfill\cr
}\vs{8pt}\hfill\cr\matrix{
\SA_{a'}\!:\!\!\!\!\hfill&
L_\mu x_\nu=(\nu+\mu)x_{\mu+\nu},\,\nu\ne0,\hfill&
L_\mu v_0=\mu(\mu+a')v_\mu,\hfill&
L_\mu y_\eta=(\eta+\OVER{\mu}2)y_{\mu+\eta}
\vs{4pt}\hfill\cr&
G_\l x_\nu=y_{\l+\nu},\,\nu\ne0,\hfill&
G_\l x_0=(2\l+a')y_\l,\hfill&
G_\l y_\eta=(\eta+\l)x_{\l+\eta}
\hfill\cr
}\vs{8pt}\hfill\cr\matrix{
\SB_{a'}\!:\hfill&
L_\mu x_\eta=(\eta+\OVER{\mu}2)x_{\mu+\eta},\hfill&
L_\mu y_\nu=\nu y_{\mu+\nu},\,\nu\ne-\mu,\hfill&
L_\mu y_{-\mu}=-\mu(\mu+a')y_0,
\vs{4pt}\hfill\cr&
G_\l x_\eta=y_{\l+\eta},\,\eta\ne-\l,\hfill&
G_\l x_{-\l}=(2\l+a')y_0,\hfill&
G_\l y_\nu=\nu x_{\l+\nu},
\hfill\cr}
\hfill\cr}
\eqno\matrix{(1.2{\rm a})\cr\vs{2pt}\cr\cr
(1.2{\rm b})\cr\vs{2pt}\cr\cr(1.2{\rm c})\cr}$$
for $\mu,\nu\in M,\,\l,\eta\in s+M$.
In this paper, we shall prove the following theorem which was conjectured
by the author in [3].
\par\ni
{\bf Theorem 1.1}. A uniformly bounded simple $\SVir$-module is
a module of the intermediate series.
\par
A $\SVir$-module $V$ is called a {\it generalized highest module} if
it is generated by a weight vector $v$ such that
there exists a \Z-basis $B=\{d_1,...,d_n\}$
of $M$ such that
$$L_\mu v=G_\l v=0,\ \forall\,
\mu=\sum_{i=1}^n m_id_i\in M\bs\{0\},\
\l=\sum_{i=1}^n n_id_i\in s+M\bs\{0\},
\eqno(1.3)$$
with all coefficients $m_i\in\Z_+,\,n_i\in\OVER12\Z_+.$
The vector $v$ is called a {\it generalized highest weight vector}.
Observe that a lowest weight module over the super-Virasoro algebras
(i.e., when $n=1$) is a generalized highest weight module by this definition.
In general, if a module $V$ has a weight
space decomposition with all the dimensions of weight spaces being finite,
then we shall call $V$ a {\it weight module}.
The following theorem generalizes Kac's conjecture in some sense.
\par\ni
{\bf Theorem 1.2}. A weight $\SVir$-module is either
a module of the intermediate series or a generalized highest weight module.
\par
Thus the problem of the classification of simple weight modules
is reduced to the problem of the classification of the generalized
highest weight modules.
\par\
\vs{-3pt}
\par
\cl{\bf\S2. Uniformly bounded modules}
\par
The goal of this section is to prove Theorem 1.1. We shall use induction
on the rank of $M$, the results of [11],
the structure of modules of the intermediate series
defined in (1.2), and the defining relations (1.1)
to obtain our result.
\par
Thus suppose that $V$ is a uniformly bounded \SVir-module.
Obviously $c$ must act trivially on such a module [3]. Thus we shall
ALWAYS omit $c$ in this section.
Decompose $V=V^{(0)}\oplus V^{(1)}$ such that
$$
\mbox{\bf\sl SVir}{\ssc\,}_iV^{(j)}
\subset V^{(i+j)},\
V^{(i)}=\bigoplus_{\mu\in\sC}V^{(i)}_\mu,
\mbox{ where }
V^{(i)}_\mu=\{v\in V^{(i)}\,|\,L_0v=\mu v\},
$$
for $i,j\in\Z/2\Z$.
For $a\in\C$, define a submodule
$$
V(a)=\sum_{\mu\in M}V^{(0)}_{\mu+a}\bigoplus\sum_{\l\in s+M}V^{(1)}_{\l+a},
$$
for $a\in\C$. Then $V$ is a direct sum of different $V(a)$. Since
we assume that $V$ is simple, we have
$V=V(a)$
for some $a\in\C.$
Since $V^{(i)},i\in\Z/2\Z$ are $\SVio$-modules, by [11, Lemma 2.1], there
exist nonnegative integers $N^{(0)},N^{(1)}$ such that
$$
{\rm dim{\sc\,}}V^{(i)}_{a+is+\mu}=N^{(i)}
\mbox{ for all $\mu\in M$ with }a+is+\mu\ne0,\,i=0,1.$$
By interchanging $V^{(0)}$ with
$V^{(1)}$ if necessary, we always suppose that $N^{(0)}\le N^{(1)}$.
\par\ni
{\bf Proof of Theorem 1.1}.
The result is obvious if $N^{(1)}=0$. Suppose that $N^{(1)}\ge1$.
If $n=1$, the result follows from [10]. So assume that $n\ge 2$.
Observe that
$${\sl SVir}[M,s]\cong{\sl SVir}[M',s']\ \Lra\
\exists\,\a\in\C\mbox{ such that }M'=\a M
\mbox{ and }s'-\a s\in M'.
$$
Thus we can always suppose that $1\in M$ and
furthermore we can suppose that $1$ is a basis element of $M$
whenever necessary.
So, we can take a \Z-basis $B=\{d_1=1,\,d_2,\cdots,d_{n-1},d_n=d\}$
of $M$ such that
$$
M=M_1\bigoplus \Z d,\ \ M_1=\Z\bigoplus\sum_{i=2}^{n-1}\Z d_i,
\ \ 2s\in M_1.
\eqno(2.1)$$
Let \SVira\, be the rank $n-1$ super-Virasoro algebra generated by
$$
\{L_\mu,G_\l\,|\,\mu\in M_1,\l\in s+M_1\}.
$$
For $k\in\Z$, let
$$
V[k]=\sum_{\mu\in M_1}V^{(0)}_{a+kd+\mu}\bigoplus
\sum_{\l\in s+M_1}V^{(1)}_{a+kd+\l},
$$
which is a uniformly bounded \SVira-module.
Since $\Z d\cap M_1=\{0\}$, there exists at most one
$k\in\Z$ such that
$a+kd\in M_1$. Fix a $k_0$ such that
$$
a+kd\not\in M_1,\ \ \forall\,k\ge k_0,
\eqno(2.2)$$
i.e., $a+kd+\nu\ne0$ for all $k\ge k_0,\,\nu\in M_1$.
First we take $k=k_0$. By inductive assumption, we can suppose that
Theorem 1.1 holds for the rank $n-1$ super-Virasoro algebra \SVira.
Thus there exists a simple \SVira-submodule $V[k,1]$ of $V[k]$ such
that $V[k,1]$ has the
form $\SA_{a+kd,b}$ in (1.2) (by (2.2), $\SA_{a+kd,b}$ is simple),
thus, there exists a basis
$$\{x_\mu\,|\,\mu\in M_1\}\cup
\{y_\l\,|\,\l\in s+M_1\},
\eqno(2.3)$$
of $V[k,1]$,
such that
$$
\mbox{\it(1.2a) holds for
all $\mu,\nu\in M_1,\,\l,\eta\in s+M_1$ with $a$ replaced by
$a+kd$.}
\eqno(2.4)$$
Furthermore, among all simple submodules of $V[k]$, we choose
$V[k,1]$ to be one such that the real part ${\rm Re}(b)$ of $b$
is minimum.
Take a composition series of the \SVira-module $V[k+1]$:
$$
0=V[k+1,0]\subset V[k+1,1]\subset\cdots\subset V[k+1,N]=
V[k+1],
\eqno(2.5)$$
for some $N$. From this, we can observe that $N=N^{(0)}=N^{(1)}$.
We shall choose the composition series (2.5) such that
$$
V[k+1,i]/V[k+1,i-1]\cong \SA_{a+(k+1)d,b'},
\eqno(2.6)$$
for some $b'\in\C,\ i=1,\cdots N,$ where $b'$ may depend on $i$.
Then, we have
$$
L_d V[k,1]^{(0)}\subset L_d V[k,1]\subset L_d V[k]\subset V[k+1],
\eqno(2.7)$$
where $V[k,1]^{(0)}$ is the ``even'' part of $V[k,1]$.
\par
We shall divide the proof into some claims.
\par
{\bf Claim 1}. $L_d V[k,1]^{(0)}\ne0.$
\par
Suppose $L_d V[k,1]^{(0)}=0.$ Using the fact
that
$$
L_{\mu+md}=
\prod_{i=-1}^{m-2}(\mu+id)^{-1}
({\rm ad\sc\,}L_d)^mL_\mu,
\ \
G_{s+\mu+md}=
\prod_{i=-1}^{m-2}(s+\mu+\OVER{i}2d)^{-1}
({\rm ad\sc\,}L_d)^m G_{s+\mu},
$$
we obtain
$$
L_{\mu+md}V[k,1]^{(0)}=G_{s+\mu+md}V[k,1]^{(0)}=0,\ \
\forall\,\mu\in M_1,\,m>0.
$$
Let ${\cal L}_1$ and ${\cal L}_2$
be the Lie super-subalgebras of $\SVir$
generated by
$$
\SVira\cup\{L_{-d}\}\mbox{ \ and \ }
\{L_{\mu+md},G_{s+\mu+md}\,|\,\mu\in M_1,\, m\in\Z_+\bs\{0\}\},
$$
respectively.
Then as space, we have $\SVir={\cal L}_1\oplus {\cal L}_2$.
By decomposing the universal enveloping algebra
$U(\SVir)=U({\cal L}_1)U({\cal L}_2)$, we obtain that
$$
V=U({\cal L}_1)U({\cal L}_2)V[k,1]^{(0)}=U({\cal L}_1)V[k,1]^{(0)}.
$$
{}From this,
we obtain that $V^{(i)}_{\ol\mu+is+md}=0$
for all $\mu\in M_1,\,m>0,\,i=0,1$. Thus $N=0$, contradicting the assumption
that $N=N^{(1)}\ge1$. Hence Claim 1 follows.
\par
Now from (2.5), (2.7), we can take $K\ge0$ to be the integer such that
$$L_d V[k,1]\not\subset V[k+1,K],\mbox{ but } L_d V[k,1]\subset
V[k+1,K+1].
\eqno(2.8)$$
Furthermore, we can choose (2.5) to be a composition series such that
$K$, defined by (2.8), is minimal possible among all composition series of
$V[k+1]$. Then, we can take
$$
\{x'_\mu\in V[k+1,K+1]^{(0)}\,|\,\mu\in M_1\},
\eqno(2.9)$$
such that
$$
(V[k+1,K+1]/V[k+1,K])^{(0)}={\rm span}
\{x'_\mu+V[k+1,K]^{(0)}\,|\,\mu\in M_1\},
$$
and such that there exist some $b'\in\C$ and some $a_\nu\in\C$ with
$$
L_\mu x'_\nu\equiv(a+(k+1)d+\nu+\mu b')x'_{\mu+\nu}, \ \
L_d x_\nu\equiv a_\nu x'_\nu\ \,({\rm mod\,} V[k+1,K]^{(0)}),
\eqno(2.10)$$
for all $\mu,\nu\in M_1.$
\par
Now [11, Lemmas 2.3-6] shows that $b'=b$ and that by rescaling
$x'_\nu$ if necessary, we can suppose that
$a_\nu=a+kd+\nu+bd$, and furthermore $K=0$, i.e.,
$V[k+1,K]=0$. Thus (2.10) becomes
$$
L_\mu x'_\nu=(a+(k+1)d+\nu+\mu b')x'_{\mu+\nu}, \ \
L_d x_\nu=(a+kd+\nu+bd) x'_\nu.
\eqno(2.11)$$
Since $V[k+1,1]$ is a $\SVira$-module of the intermediate series,
there exists basis $\{y'_\l\,|\,\l\in s+M_1\}$ of $V[k+1,1]^{(1)}$ such that
$$
\matrix{
\mbox{\it(1.2a) holds for all $\mu,\nu\in M_1,\,\l,\eta\in s+M_1$ with the symbols
$x,y$ replaced}
\vs{4pt}\hfill\cr
\mbox{\it by symbols $x',y'$ respectively,
and with $a$ replaced by $a+(k+1)d$.}
\hfill\cr}
\eqno(2.12)$$
By (2.8), $L_dV[k,1]^{(1)}\subset V[k+1,1]^{(1)}$.
Thus we can suppose
$$
L_d y_\l=c_\l y'_\l\mbox{ for some }c_\l\in\C,\ \l\in s+M_1.
\eqno(2.13)$$
Then by [11, Lemma 2.3], there exists nonzero scalar $c\in\C$ such that
$$
c_\l=c(a+kd+\l+d(b-\OVER12)),\ \ \forall\,\l\in s+M_1.
$$
\par
{\bf Claim 2}. One can choose a suitable simple
$\SVira$-submodule $V[k,1]$ of $V[k]$
such that $L_{-d}L_d V[k,1]\subset V[k,1]$.
\par
Let $W[k]$ be the direct sum of
ALL simple \SVira-submodules $V[k,1]$ of $V[k]$ of the same type:
$V[k,1]\cong \SA_{a+kd,b}$ and let
$W[k+1]$ be the direct sum of
ALL simple \SVira-submodules $V[k+1,1]$ of $V[k+1]$ of the same type:
$V[k+1,1]\cong \SA_{a+(k+1)d,b}$.
Then the above discussion shows that
$L_d W[k]\subset W[k+1]$, and so we also have
$L_{-d} W[k+1]\subset W[k]$ by replacing $d$ by $-d$ in the above arguments.
By Claim 1, $W[k],W[k+1]$ must have the same number, say, $r$,
of composition factors.
Choose a suitable basis
$\{x^{(j)}_\nu,y^{(j)}_{s+\nu}\,|\,\nu\in M_1, j=1,...,r\}$ of $W[k]$
such that $x^{(1)}_0$ is an eigenvector of $L_{-d}L_d$ and
$$
\matrix{
L_\mu X_\nu=(a+kd+\nu+\mu b)X_{\mu+\nu},\hfill&
L_\mu Y_\eta=(a+kd+\eta+\mu(b-\OVER12))Y_{\mu+\eta}
\vs{4pt}\hfill\cr
G_\l X_\nu=Y_{\l+\nu},\hfill&
G_\l Y_\eta=(a+kd+\eta+2\l(b-\OVER12))X_{\l+\eta},
\hfill\cr}
\eqno(2.14)$$
where
$$
X_\nu=\pmatrix{x_\nu^{(1)}\cr\vdots\cr x_\nu^{(r)}\cr},\ \ \
Y_\eta=\pmatrix{y_\eta^{(1)}\cr\vdots\cr y_\eta^{(r)}\cr},\ \ \
\forall\,\nu\in M_1,\ \eta\in s+M_1.
\eqno(2.15)$$
Then Claim 1, (2.11) and (2.13) allow us to choose a basis
$\{x'^{(j)}_\nu,y'_{s+\nu}\,|\,\nu\in M_1,j=1,...,r\}$ of $W[k+1]$ such
that
$$
L_d X_\nu=(a+kd+\nu+bd)X'_\nu,\ \ \
L_d Y_\eta=c(a+kd+\eta+d(b-\OVER12))Y'_\eta,
$$
where $X'_\nu,Y'_\eta$ are the similar notations as in (2.15),
and such that
(2.14) holds for $k$ replaced by $k+1$ and $X,Y$ replaced by $X',Y'$.
Similarly, with $d$ replaced by $-d$ in the above arguments, we have
$$
L_{-d} X'_\nu=(a+(k+1)d+\nu-bd)A X_\nu,\ \ \
L_{-d} Y'_\eta=c'(a+(k+1)d+\eta-d(b-\OVER12))AY_\eta,
$$
where $A$ is some constant $r\times r$ matrix and $c'\in\C\bs\{0\}$.
Thus we have
$$
L_{-d}L_d X_\nu=(a+kd+\nu+bd)(a+(k+1)d+\nu-bd)AX_\nu.
\eqno(2.16)$$
Since $x^{(1)}_0$ is an eigenvector of $L_{-d}L_d$, on the first row of
$A$, by (2.16),
there is only one possible nonzero entry in the first position.
If we choose
$$
V[k,1]={\rm span}_\sC\{x^{(1)}_\nu,y^{(1)}_{s+\nu}\,|\,\nu\in M_1\},
$$
then we see that $L_{-d}L_d V[k,1]\subset V[k,1]$.
This proves our Claim 2.
\par
Now we can complete the proof of Theorem 1.1 as follows.
Starting from $k=k_0$ and choosing $V[k_0,1]$ as in Claim 2,
by induction on $k>k_0$, we can choose $V[k,1]$ to be the simple
\SVira-submodule generated by $L_d V[k-1,1]$.
As in the proof of Claim 1, using
$L_{-d}L_d V[k_0,1]\subset V[k_0,1]$, by induction on $m$,
we have
$$
L_{\nu+md}V[k,1]\subset V[k+m,1],\
G_{\l+md}V[k,1]\subset V[k+m,1],
\eqno(2.17)$$
for all $\nu\in M_1,\,\l\in s+M_1,\,m,k\in\Z$
such that $k\ge k_0,\,k+m\ge k_0$.
Since $V$ is a simple \SVir-module, it must be generated by
$V[k_0,1]$, thus by (2.17),
$$
{\rm dim\sc\,}V^{(0)}_{a+\nu+kd}=
{\rm dim\sc\,}V^{(1)}_{a+\l+kd}=1,
\ \ \forall\,\nu\in M_1,\,\l\in s+M_1,\,k\ge k_0.
$$
This shows that $N=1$, thus
${\rm dim\sc\,}V^{(0)}_{a+\nu}=1=
{\rm dim\sc\,}V^{(1)}_{a+\l}$
for all $\nu\in M,\,\l\in s+M$ with
$a+\nu\ne0\ne a+\l$, from this one sees that
we must also have
${\rm dim\sc\,}V^{(0)}_{a+\nu}\le1,\,
{\rm dim\sc\,}V^{(1)}_{a+\l}\le1$ if $a+\nu=0=a+\l$.
Thus $V$ must be a module of the intermediate series.
\par\
\vs{-3pt}
\par
\cl{\bf\S3. Generalized highest weight modules}
\vs{-1pt}
\par
The aim of this section is to prove Theorem 1.2. This can be done
by three lemmas below. We suppose $n\ge2$ as for $n=1$, the result
follows from [10].
\par\ni{\bf Lemma 3.1}.
Let $V$ be a simple \SVir-module and let $B=\{d_1,...,d_n\}$ be
any \Z-basis of $M$.
If $0\ne v\in V$ is a weight vector and $k\in\Z_+$ such that
$L_\mu v=0=G_\l v$ for all
$$
\matrix{
\dis\mu\in M_{\ssc B}^{(k)}=
\{\mu\in M\,|\,\mu=\sum_{i=1}^n m_id_i,
\mbox{ all }k\le m_i\in\Z_+\},
\vs{4pt}\hfill\cr
\dis\l\in M'^{(k)}_{\ssc B}=\{\l\in s+M\,|\,
\l=\sum_{i=1}^n n_id_i,\mbox{ all }k\le m_i\in\OVER12\Z_+\},
\hfill\cr}
$$
then $V$ is a generalized highest weight module.
\par\ni{\bf Proof}. Take
$$
d'_i=\sum_{j=1}^i(k+i-j+1)d_j+k\sum_{j=i+1}^nd_j,\,i=1,...,n.
$$
It is straightforward to check that the determinant of coefficients of
$\{d_1,...,d_n\}$ is 1, thus
$B'=\{d'_1,...,d'_n\}$ is also a \Z-basis of $M$.
Now for any
$$
0\ne\mu'=\sum_{i=1}^n m'_i d'_i\in M_{\ssc{B'}}^{(0)}
\mbox{ (so all $m'_i\ge0$, at least one $m'_i>0$),}
$$
we have $\mu'=\sum_{i=1}^n m_id_i$ with all coefficients
$$
m_i=k\sum_{j=1}^{i-1}m'_j+
\sum_{j=i}^n(k+j+1-i)m'_j\ge k,
$$
i.e., $\mu'\in M_{\ssc B}^{(k)}$,
thus $L_{\mu'\ssc\,}v=0$.
Similarly, for any $0\ne\l'\in M'^{(0)}_{\ssc{B'}}$, we have
$G_{\l'\ssc\,} v=0$. Hence by definition (cf. (1.3)$\ssc\,$),
$v$ is a generalized highest weight vector and $V$ is a
generalized highest weight module.
\qed\par\ni
{\bf Lemma 3.2}.
Fix a \Z-basis $B=\{d_1,...,d_n\}$ of $M$.
Let
$$
\matrix{
\dis A=\{\mu\in M\,|\,\mu=\sum_{i=1}^nm_id_i,\mbox{ all }m_i\in\{-1,0,1\}\},
\vs{4pt}\hfill\cr
\dis A'=\{\l\in s+M\,|\,\l=\sum_{i=1}^nn_id_i,\mbox{ all }
n_i\in\{-1,-\OVER12,0,\OVER12,1\}\},
\hfill\cr}
$$
be two finite sets.
For any $\mu=\sum_{i=1}^nm_id_i\in M$,
let
$$
A_\mu=\{\nu\in M\,|\,\nu-\mu\in A\},\ \
A'_\mu=\{\l\in s+M\,|\,\l-\mu\in A'\},
$$
\def\LL{\mbox{$\cal L$}}and
let $\LL$ be the Lie super-subalgebra of $\SVir$ generated by
$\{L_\nu\,|\,\nu\in A_\mu\}\cup\{G_\l\,|\,\l\in A'_\mu\}$.
Then there exists a \Z-basis $B'=\{d'_1,...,d'_n\}$ of $M$
such that
$$
\LL\supset \SVir_{\ssc{B'}}^{(k)}=
\{L_\nu,\,G_\l\,|\,\nu\in M_{\ssc{B'}}^{(k)},\,\l\in
 M'^{(k)}_{\ssc{B'}}\}\mbox{ for some }
k\in\Z_+.
$$
\par\ni{\bf Proof}.
If necessary, by replacing $d_i$ by $-d_i$ (this does not change the
sets $A,A'$, neither do $\LL$), we can suppose $m_i\ge0,$ $i=1,...,n$.
First suppose that $m_1\ne0,m_2\ne0$. Take
$$
\matrix{
d'_1=\hfill&\!\!\!m_2\mu\!+\!d_1\hfill&\!\!\!=\hfill&\!\!\!
\dis(\mu+d_1)+\sum_{i=1}^{m_2-1}\mu=\hfill&\!\!\!(m_1m_2+1)d_1+m^2_2d_2+
\dis\sum_{i=3}^nm_2m_id_i,
\hfill\cr}
$$
$$
\matrix{
d'_2=\hfill&\!\!\!m_1\mu\!-\!d_2\hfill&\!\!\!=\hfill&\!\!\!
\dis(\mu-d_2)+\sum_{i=1}^{m_1-1}\mu=\hfill&\!\!\!m^2_1d_1+(m_1m_2-1)d_2+
\dis\sum_{i=3}^nm_1m_id_i,\vs{5pt}\hfill\cr
d'_i=\hfill&\!\!\!d'_1+d_i\hfill&\!\!\!=\hfill&\!\!\!
\dis(\mu\!+\!d_1\!+\!d_i)+\sum_{i=1}^{m_2-1}\mu,\!\!\!\!\hfill&i=3,...,n.\hfill
\cr}$$
We can easily solve $d_i,i\!=\!1,...,n,$ as an
integral linear combination of $B'\!=\!\{d'_1,...,d'_n\}$, so $B'$ is a
\Z-basis of $M$. From the second equality of $d'_i$, we see
$L_{d'_i}\!\in\!\LL,\,i\!=\!1,...,n$. For example,
$$
[L_\mu,...,[L_\mu,L_{\mu+d_1}]...]\ (m_2-1
\mbox{ copies of }L_\mu)= aL_{d'_1},
\mbox{ where }a=\prod_{i=0}^{m_2-2}(i\mu+d_1)\ne0,
$$
so that the left-hand side is in $\LL$, thus $L_{d'_1}\in \LL$ and
so
$\{L_\nu\,|\,\nu\in M_{\ssc{B'}}^{(0)}\}\subset\LL$. Similarly,
$\{G_\l\,|\,\l\in M'^{(0)}_{\ssc{B'}}\}\subset\LL$.
Next suppose that some $m_i=0$, $1\le i\le n$, say, $m_1=0$. Take
$$
d'_1=\mu+d_1,\ \ d'_i=d'_1+d_i=\mu+d_1+d_i, \ i=2,...,n.
$$
Then as above $B'=\{d'_1,...,d'_n\}$ is a \Z-basis of $M$ and
$\HVir_{_{B'}}^{(0)}\subset\LL$.
\qed\par
\ni{\bf Lemma 3.3}.
A simple weight \SVir-module $V$ with dimensions of weight spaces
not uniformly bounded is a generalized highest weight module.
\par\ni
{\bf Proof}.
Suppose conversely that $V$ is not a
generalized highest weight module.
For any $\mu\in M$, let $A,\,A',\,A_\mu,\,A'_\mu,\,\LL$
and $B'$ be as in Lemma 3.2. By Lemmas 3.1 and 3.2, for any
$0\ne v\in V$, we have $\SVir_{\ssc{B'}}^{(k)}v\ne0$, thus,
$L_\nu v\ne0$ for some $\nu\in A_\mu$ or
$G_\l v\ne0$ for some $\l\in A'_\mu$, i.e.,
$$
\bigcap_{\nu\in A_\mu}{\rm ker}(L_\nu|_V)
\bigcap\bigcap_{\l\in A'_\mu}{\rm ker}(G_\l|_V)=0.
$$
In particular,
$$
(\bigoplus_{\nu\in A_\mu}L_\nu\bigoplus\bigoplus_{\l\in A'_\mu}G\l)
|_{V_{a-\mu}}:\
V_{a-\mu}\rar\bigoplus_{\nu\in A}V_{a+\nu}
\bigoplus\bigoplus_{\l\in A'}V_{a+\l},
$$
is an injection.
Therefore, ${\rm dim\sc\,}V_{a-\mu}\le N$, where
$$
N=
\sum_{\nu\in A}{\rm dim\sc\,}V_{a+\nu}
+\sum_{\l\in A'}{\rm dim\sc\,}V_{a+\l},
$$
is a fixed finite integer since both $A$ and $A'$ are finite sets.
As $\mu\in M$ is arbitrary,
this proves that $V$ is uniformly bounded, a
contradiction.
\qed\par
Now Theorem 1.2 follows from Lemma 3.3.
\par\
\vs{-3pt}
\par
\cl{\bf References}
\par\small
\parskip .03 truein
\baselineskip 4pt
\lineskip 4pt
\ni\hi3ex\ha1
1 \ Gelfand,~I.~M., Fuchs,~D.~B.: Cohomologies of the Lie algebra of vector
 fields on the circle, {\sl Funct.~Anal.~Appl.}, {\bf2}\,(1968),~92-39
 (English translation 114-126).
\par\ni\hi3ex\ha1
2 \ Patera,~J., Zassenhaus,~H.: The higher rank Virasoro algebras,
 {\sl Comm.~Math.~Phys.}, {\bf 136}\,(1991), 1-14.
\par\ni\hi3ex\ha1
3 \ Su,~Y.: Harish-Chandra modules of the intermediate series over
 the high rank Virasoro algebras and high rank super-Virasoro
 algebras, {\sl J.~Math.~Phys.}, {\bf 35}\,(1994), 2013-2023.
\par\ni\hi3ex\ha1
4 \ Su,~Y., Zhao,~K.: Generalized Virasoro and super-Virasoro algebras
 and modules of the intermediate series, to appear.
\par\ni\hi3ex\ha1
5 \ Kac,~V.~G.: Some problems of infinite-dimensional Lie algebras and their
 representations, in: {\sl Lect. Notes in Math.}, {\bf933}, 117-126.
 Berlin, Heidelberg, New York: Springer, 1982.
\par\ni\hi3ex\ha1
6 \ Chari,~V., Pressley,~A.: Unitary representations of the Virasoro
 algebra and a conjecture of Kac, {\sl Compositio Math.},
 {\bf 67}\,(1988), 315-342.
\par\ni\hi3ex\ha1
7 \ Martin,~C., Piard,~A.: Indecomposable modules over the Virasoro
 Lie algebra and a conjecture of V.~Kac, {\sl Comm.~Math.~Phys.},
 {\bf 137}\,(1991), 109-132.
\par\ni\hi3ex\ha1
8 \ Su,~Y.: A classification of indecomposable $sl_2(\C)$-modules
 and a conjecture of Kac on Irreducible modules over the Virasoro
 algebra, {\sl J.~Alg.}, {\bf 161}\,(1993), 33-46.
\par\ni\hi3ex\ha1
9 \ Mathieu,~O.: Classification of Harish-Chandra modules over the
 Virasoro Lie algebra, {\sl Invent.~Math.}, {\bf 107}\,(1992), 225-234.
\par\ni\hi3ex\ha1
10 Su,~Y.: Classification of Harish-Chandra modules over the
 super-Virasoro algebras, {\sl Comm. Alg.}, {\bf 23}\,(1995), 3653-3675.
\par\ni\hi3ex\ha1
11 Su,~Y.: Simple modules over the High Rank Virasoro algebras,''
 {\sl Comm.~Alg.}, in press.
\par\ni\hi3ex\ha1
12 Su,~Y: On indecomposable modules over the
 Virasoro algebra, {\sl Sciences in China A}, in press.
%
\end{document}